\documentclass[12pt,a4paper,reqno]{amsart}
\usepackage[latin1]{inputenc}
\usepackage{amsmath}
\usepackage{amsfonts}
\usepackage{amssymb}
\usepackage{amsmath,bbm,amssymb,amsxtra}
\usepackage{enumerate}
\usepackage{caption}
\allowdisplaybreaks
\usepackage{epsfig}
\usepackage{ae}
    \def\Xint#1{\mathchoice
    {\XXint\displaystyle\textstyle{#1}}%
    {\XXint\textstyle\scriptstyle{#1}}%
    {\XXint\scriptstyle\scriptscriptstyle{#1}}%
    {\XXint\scriptscriptstyle\scriptscriptstyle{#1}}%
    \!\int}

    \def\XXint#1#2#3{{\setbox0=\hbox{$#1{#2#3}{\int}$ }
    \vcenter{\hbox{$#2#3$ }}\kern-.6\wd0}}
    \def\dashint{\Xint-}

\theoremstyle{definition}

\newtheorem{lemma}{Lemma}[section]

\newtheorem{proposition}[lemma]{Proposition}

\newtheorem{theorem}[lemma]{Theorem}

\newtheorem{corollary}[lemma]{Corollary}

\newtheorem{remark}[lemma]{Remark}

\newtheorem*{acknowledgements}{Acknowledgements}
\newtheorem{definition}[lemma]{Definition}

\newcommand{\prop}[1]{\begin{proposition}\label{#1}
\sl }
\newcommand{\eprop}{\end{proposition}}
\newcommand{\thm}[1]{\begin{theorem}\label{#1}
\sl }
\newcommand{\ethm}{\end{theorem}}
\newcommand{\cor}[1]{\begin{corollary}\label{#1}
\sl }
\newcommand{\ecor}{\end{corollary}}

\newcommand{\lem}[1]{\begin{lemma}\label{#1}
\sl }
\newcommand{\elem}{\end{lemma}}
\newcommand{\defin}[1]{\begin{definition}\label{#1}
\sl }
\newcommand{\edefin}{\end{definition}}
\newcommand{\beqno}{\begin{eqnarray*}}
\newcommand{\eeqno}{\end{eqnarray*}}
\newcommand{\beqla}[1] {\begin {eqnarray}\label{#1}}
\def\eeq {\end {eqnarray}}
\newcommand{\beq}{\begin {eqnarray}}
\newcommand{\real}{{\mathbb R}}

\newcommand{\integer}{{\mathbb Z}}

\newcommand{\nanu}{{\mathbb N}}

\newcommand{\complex}{{\mathbb C}}

\newcommand{\trace}{{\mathbb T\hskip-4pt{\mathbb R}} }

\newcommand{\integral}{{\mathbb I}}

\newcommand{\Lip}{{\rm Lip}\,}
\newcommand{\D}{\mathbb{D}}
\newcommand{\hlmax}{\mathcal{M}}

\newcommand{\df}{\dot{\mathcal F}}
\newcommand{\dfp}{\dot{\mathcal F}^{s}_{p,q}(Z)}

\newcommand{\db}{\dot{\mathcal B}}
\newcommand{\dbp}{\dot{\mathcal B}^{s}_{p,q}(Z)}

\newcommand{\dn}{\dot{N}}
\newcommand{\dnp}{\dot{N}^s_{p,q}(Z)}

\newcommand{\J}{{\mathcal J}}
\newcommand{\I}{{\mathcal I}}

\newcommand{\jp}{\J^{s}_{p,q}(X)}

\newcommand{\ip}{\I^{s}_{p,q}(X)}
\newcommand{\ipe}{\I^{s}_{p,q}(E)}

\newcommand{\izero}{\mathring{\I}^s_{p,q}(X)}
\newcommand{\differencespace}{\mathring{{\mathcal D}}\jp}
\newcommand{\differencespacei}{\mathring{{\mathcal D}}\ip}
\newcommand{\locint}{L^1_{\rm{loc}}}

\DeclareMathOperator*{\esssup}{ess\,sup}

\title[Besov spaces]{Besov spaces via hyperbolic fillings}

\author[Soto]{Tom\'as Soto}
\address{Department of Mathematics and Statistics, University of Helsinki, PO~Box~68, FI-00014 Helsinki, Finland}
\email{tomas.soto@helsinki.fi}
\thanks{The author was supported by the Finnish CoE in Analysis and Dynamics Research}



\begin{document}

\maketitle

\begin{abstract}
We establish a new characterization of the homogeneous Besov spaces $\dbp$ with smoothness $s \in (0,1)$ in the setting of doubling metric measure spaces $(Z,d,\mu)$. The characterization is given in terms of a \emph{hyperbolic filling} of the metric space $(Z,d)$, a construction which has previously appeared in the context of other function spaces in \cite{BP,BS,BSS}. We use the characterization to obtain results concerning the density of Lipschitz functions in the spaces $\dbp$ and a general complex interpolation formula in the smoothness range $0 < s < 1$.
\end{abstract}

\section{Introduction}\label{se:introduction}

In \cite{BSS}, we formulated a theory of the family of function spaces $\df^s_{p,q}$, known as the Triebel-Lizorkin spaces, in the setting of homogeneous metric measure spaces using a construction called the \emph{hyperbolic filling} of the metric space. In the context of function spaces, this construction has previously appeared in e.g.~\cite{BP,BS}. The purpose of this note is to complement \cite{BSS}; we formulate a basic theory of the \emph{Besov spaces} $\db^s_{p,q}$ in the smoothness range $0 < s < 1$ via a hyperbolic filling of the underlying metric space. A similar family of function spaces in the parameter range $p = q$ appeared in \cite{BP}.

Let us start by explaining the assumptions on the metric measure spaces under consideration. Throughout the paper, $Z := (Z,d,\mu)$ is assumed to be a (bounded or unbounded) metric measure space such that the measure $\mu$ is Borel regular and all balls (with respect to the metric $d$) have positive and finite $\mu$-measure. $\mu$ is also assumed to have the following \emph{doubling property}: there exists a finite constant $c \geq 1$ such that $\mu\big(B(\xi,2r)\big) \leq c B\mu\big(B(\xi,r) \big)$ for all $\xi \in Z$ and $r > 0$. From this assumption it follows that there exist constants $C \geq 1$ and $Q > 0$ such that
\beqla{eq:doubling}
  \mu\big( B(\xi,\lambda r) \big) \leq C \lambda^Q \mu\big( B(\xi,r)\big)
\eeq
for all $\xi \in Z$, $r > 0$ and $\lambda \geq 1$. $Q$ is in some sense an upper bound for the dimension of $Z$, and it will be fixed from now on.

We shall next explain the construction of the hyperbolic filling of $Z$. For all $n \in \integer$, let $(\xi_x)_{x \in X_n}$, where $X_n$ is a suitable index set, be a maximal set of points of $Z$ with pairwise distances at least $2^{-n-1}$. For $x \in X_n$, we write $|x|$ for $n$, the \emph{level} of $x$. It now follows easily from the doubling property of the underlying space that for all $n \in \integer$, the balls $B(x) := B(\xi_x,2^{-n})$ with $x \in X_n$ have uniformly bounded overlap in $Z$, that the balls $B(\xi_x,2^{-n-1})$ with $x \in X_n$ still cover $Z$ and that the balls $B(\xi_x,2^{-n-2})$ with $x \in X_n$ are pairwise disjoint.

Write $X := \sqcup_{n \in \integer} X_n$, and write $(X,E)$ for graph with the property that two distinct vertices $x$, $x' \in X$ are joined by an edge in $E$ if and only if $B(x)\cap B(x') \neq \emptyset$ and $||x|-|x'|| \leq 1$; in this case we also write $x \sim x'$. Under some mild additional assumptions on $Z$, the natural path metric on the graph $(X,E)$ makes it hyperbolic in the sense of Gromov so that its boundary at infinity coincides with $Z$ (see \cite{BP}), which is why it is called a hyperbolic filling of $Z$.

Now for a complex-valued sequence $u$ on $X$, the \emph{discrete derivative} $du$ could be defined as a scalar-valued sequence on $E$ or, essentially in an equivalent way, as a vector-valued sequence on $X$. To simplify our notation, we shall for now take the latter approach, and postpone the former approach until Section \ref{se:complex-interpolation} below. More precisely, the doubling assumption implies that the graph $(X,E)$ has bounded valency, i.e.~there exists $\Delta(X) \in \nanu$ such that every point of $X$ is joined to at most $\Delta(X)$ distinct points of $X$ by an edge in $E$. Then $du$ is defined as an $\complex^{\Delta(X)}$-valued sequence on $X$ by setting
\[
  du (x) := \big(u(y_1) - u(x) , u(y_2) - u(x) , \cdots , u(y_{\deg(x)}) - u(x), 0 , \cdots \big)
\]
for all $x \in X$, where $y_1$, $\cdots$, $y_{\deg(x)}$ are the neighbors of $x$ in some order, and the vector above is augmented by zeroes in case $\deg(x) < \Delta(X)$. In particular,
\[
  |du(x)| =  \Big(\sum_{x' \sim x} |u(x') - u(x)|^2 \Big)^{1/2}
\]
for all $x \in X$.

We write $\locint(Z)$ for the vector space of complex-valued measurable functions on $Z$ that are integrable on bounded subsets of $Z$. For $f \in \locint(Z)$, the \emph{Poisson extension} $Pf\colon X\to \complex$ is defined by
\[
  Pf(x) := \dashint_{B(x)} f\, d\mu := \frac{1}{\mu\big(B(x)\big)} \int_{B(x)} f\, d\mu.
\]
Then, our Besov space $\dbp$ is defined as the space of functions $f \in \locint(Z)$ such that
\[
  \| f \|_{\dbp} := \big\|\, |d(Pf)|  \,\big\|_{s,p,q}
\]
is finite, where $\| \cdot\|_{s,p,q}$ is a certain sequence norm on $X$. An in-depth discussion of the motivation behind a definition like this can be found in the introduction of \cite{BSS}. For now, let us simply note that in the smoothness range $0 < s < 1$, this definition coincides with several other definitions that have appeared in previous literature, including the standard Fourier-analytical definition in the Euclidean setting; see Proposition \ref{pr:hajlasz} below.

Before going to the precise definitions of our spaces, let us introduce some notation. If $B$ is a ball in $Z$ with a distinguished center point and a radius, we write $\lambda B$ (where $\lambda > 0$) for the ball with the same center point and radius $\lambda$ times the original radius. The notations $\locint(Z)$ and $\dashint_E f d\mu$, where $E$ is a measurable subset of $Z$ and $f$ is a complex-valued function on $Z$, will be used with the same meanings as above. For any two non-negative functions $f$ and $g$ defined on the same set, the notation $f \lesssim g$ means that $f \leq C g$, where $C$ is a finite constant usually independent of some parameters that will be obvious from the context. The notation $f \approx g$ means that $f \lesssim g$ and $g \lesssim f$.

\section{Definitions and basic properties}\label{se:basic}

The results presented here are analogs to the corresponding results for the Triebel-Lizorkin spaces $\dfp$ presented in \cite{BSS}. Let us start by giving the definition of the sequence space on which our Besov spaces will be based on.

\defin{de:besov}
Let $0 < s < \infty$ and $0 < p,\,q \leq \infty$. $\ip$ is the quasi-normed space of sequences $u\colon X \to \complex$ such that
\begin{align}
  \|u\|_{\ip} & := \bigg(\sum_{k\in\integer} 2^{ksq} \big\| \sum_{x \in X_k} |u(x)| \chi_{B(x)} \big\|_{L^p(Z)}^q\bigg)^{1/q} \notag \\
& \approx \bigg( \sum_{k\in\integer} 2^{ksq }\Big(\sum_{x \in X_k} \mu\big(B(x)\big) |u(x)|^p \Big)^{q/p} \bigg)^{1/q} \label{eq:seq-norm}
\end{align}
(obvious modifications for $q = \infty$ and/or $p = \infty$) is finite.
\edefin

\prop{pr:quasi-banach}
{\rm (i)} $\ip$ is a quasi-Banach space for all admissible parameters. When $1 < p,\,q < \infty$, it is a reflexive Banach space.

\smallskip
{\rm (ii)} Let $(c_x)_{x \in X}$ be a sequence of positive numbers such that $0 < \inf_{x \in X} c_x \leq \sup_{x \in X} c_x < \infty$. Then
\[
  u \mapsto \bigg( \sum_{k\in\integer} 2^{ksq }\Big(\sum_{x \in X_k} \mu\big(c_x B(x)\big) |u(x)|^p \Big)^{q/p} \bigg)^{1/q}
\]
(obvious modification for $q = \infty$ and/or $p = \infty$) is an equivalent quasinorm on $\ip$.
\eprop

\begin{proof}
(i) As seen from \eqref{eq:seq-norm}, the quasinorm of $\ip$ is essentially equivalent to a weighted quasinorm of the type $l^q(l^p)$, so this can be proven by a standard argument.

(ii) This is an immediate consequence of \eqref{eq:seq-norm} and the doubling property of $\mu$.
\end{proof}

Before giving the definition of the spaces $\dbp$, let us formulate several important auxiliary results. We shall first examine the boundary behavior of sequences on $X$ that essentially correspond to our Besov spaces. For this purpose, fix a collection $(\psi_x)_{x \in X}$ of non-negative Lipschitz functions such that $\psi_x$ is supported on $B(x)$ for all $x$, $(\psi_x)_{x \in X_n}$ is a partition of unity in $Z$ for all $n \in \integer$ and $\Lip \psi_x \lesssim 2^{|x|}$ for all $x$. For a sequence $u$ on $X$ we then define
\[
  T_n u := \sum_{x \in X_n} u(x) \psi_x
\]
for all $n \in \integer$. The following result is an analog to \cite[Lemma 2.3]{BSS}.

\lem{le:trace}
Suppose that $0 < s < \infty$, $Q/(Q+s) < p \leq \infty$ and $0 < q \leq \infty$. If $u$ is a sequence on $X$ such that $|du| \in \ip$, then the limit
\[
  \trace u := \lim_{n \to \infty} T_n u
\]
exists in $\locint(Z)$ and pointwise $\mu$-almost everywhere, and
\[
  \|\trace u - u(x)\|_{L^1(B(x))} \lesssim \|du\|_{\ip}
\]
for every $x \in X$, where the implied constant depends on $x$ but not on $u$.
\elem
\begin{proof}
The proof is similar to the proof of \cite[Lemma 2.3]{BSS}. Fix $x$ and let $N := |x|$. Then a simple calculation yields
\begin{align}
  & \int_{B(x)} |T_{N}u - u(x)| d\mu + \sum_{n \geq N} \int_{B(x)} |T_{n+1}u - T_n u| d\mu \label{eq:trace-conv}\\ 
&  \qquad
\lesssim \sum_{|y| \geq N,\; B(y) \cap B(x) \neq \emptyset} \mu\big(B(y)\big) |du(y)|
\leq \bigg( \sum_{|y| \geq N,\; B(y) \cap B(x) \neq \emptyset} \mu\big(B(y)\big)^r |du(y)|^r \bigg)^{1/r}, \notag
\end{align}
where $r$ is any positive number strictly less than $1$. In particular, we can take $\epsilon \in (0,s)$ so that $r := Q/(Q+\epsilon) < \min(1,p)$. Now $r = 1 - (\epsilon/Q)r$ and the doubling property yields $\mu\big(B(y)\big)^{-(\epsilon/Q)r} \lesssim \big[2^{(|y|-|x|)Q(\epsilon/Q)} \mu\big(B(x)\big)^{-\epsilon/Q}\big]^r \approx 2^{|y|\epsilon r}$ for all $y$ in the sum above (where the last implied constant depends on $x$). Since also $\epsilon < s$, we may estimate the sum by a constant times
\begin{align*}
&  \bigg( \sum_{|y| \geq N,\; B(y) \cap B(x) \neq \emptyset} 2^{|y|\epsilon r}\mu\big(B(y)\big) |du(y)|^r \bigg)^{1/r}\\
& \qquad \lesssim \bigg( \sum_{k \geq N} 2^{ksq} \Big(\sum_{y \in X_k,\; B(y)\cap B(x) \neq \emptyset} \mu\big(B(y)\big) |du(y)|^r \Big)^{q/r} \bigg)^{1/q}.
\end{align*}
Finally, since $p/r > 1$ and the balls corresponding to any level of $X$ have uniformly bounded overlap, the $k$th term in the outer sum above can be estimated using H\"older's inequality as follows:
\[
  2^{ksq} {\underbrace{\Big(\sum_{y \in X_k,\; B(y)\cap B(x) \neq \emptyset} \mu\big(B(y)\big) \Big)}_{\approx \mu(B(x)) \approx 1}}^{\frac{q}{r} \frac{p/r - 1}{p/r}} \Big(\sum_{y \in X_k,\; B(y)\cap B(x) \neq \emptyset} \mu\big(B(y)\big) |d u(y)|^p \Big)^{q/p}.
\]
All in all, the left-hand side of \eqref{eq:trace-conv} can be estimated from above by a constant times
\[
  \bigg( \sum_{k \geq N} 2^{ksq} \Big(\sum_{y \in X_k,\; B(y)\cap B(x) \neq \emptyset} \mu\big(B(y)\big) |du(y)|^p \Big)^{q/p} \bigg)^{1/q} \lesssim \|du\|_{\ip}.\qedhere
\]
\end{proof}

It turns out that if $|du| \in \ip$, then we have $P(\trace u) = u$ in some sense near the ``boundary'' of the graph $(X,E)$. This will be quantified in Theorem \ref{th:trace-commute} below, but before stating said Theorem, we will formulate the following auxiliary result, which will play an important role in many subsequent proofs.

\prop{pr:t-operator}
Let $0 < s < \infty$, $Q/(Q+s) < p \leq \infty$ and $0 < q \leq \infty$. Then the operator $T$ defined by
\[
  Tu(x) := \sum_{\substack{|y| \geq |x|\\B(y)\cap B(x) \neq \emptyset}} \frac{\mu\big (B(y)\big)}{\mu\big(B(x)\big)} u(y)
\]
is well-defined and bounded on $\ip$.

The conclusion continues to hold for the operator $u \mapsto (Tu)\circ\Psi$ whenever $\Psi \colon X\to X$ is a mapping such that $B(\Psi(x))\cap B(x) \neq \emptyset$ for all $x$ and there exists $\sigma \geq 0$ so that $||\Psi(x)| - |x|| \leq \sigma$ for all $x$.
\eprop

\begin{proof}
The necessary technical details are in essence contained in the proof of \cite[Proposition 2.4]{BSS}, so we shall only indicate how to adapt them to the situation of this proof. We will only consider the case $p < \infty$, as the case with $p = \infty$ can be handled in a similar but easier manner.

Let $u \in \ip$. For $k \in \integer$, define the function $U_k \colon Z\to [0,\infty]$ by
\[
  U_k(\xi) = \| \{2^{|x|s} |u(x)| \,:\, x \in X_k,\, B(x) \owns \xi \}\|_{\ell^p},
\]
so that the $\ip$-norm of $u$ is essentially obtained as the $l^q(L^p)$-norm of the functions $U_k$. Now taking $r$ so that $Q/(Q+s) < r < \min(1,p)$, the first part of the proof of \cite[Proposition 2.4]{BSS} establishes the following estimate for all $\xi \in Z$ and $x \in X$ such that $B(x) \owns \xi$:
\[
  2^{|x|s} T(|u|)(x) \lesssim \sum_{k \geq |x|} 2^{(|x|-k)(Q+s-Q/r)} \hlmax\big(U^r_k\big)(\xi)^{1/r},
\]
where $\hlmax$ stands for the Hardy-Littlewood maximal function and the implied constant is independent of $\xi$ and $x$. Thus if $n \in \integer$, we can use the bounded overlap of the balls corresponding to the vertices in $X_n$ together with either H\"older's inequality (if $p > 1$) or the subadditivity of the function $t \mapsto t^p$ (otherwise) to obtain
\begin{align*}
  2^{ns} \big\| \sum_{x \in X_n} T(|u|)(x) \chi_{B(x)}\big\|_{L^p(Z)} & \lesssim \big\| \sum_{k \geq n} 2^{(n-k)(Q+s-Q/r)}\hlmax\big(U^r_k\big)^{1/r} \big\|_{L^p(Z)} \\
& \lesssim \Big(\sum_{k \geq n}2^{(n-k)(Q+s-Q/r)(p\land 1)} \big\| \hlmax\big(U_k^r\big) \big\|_{L^{p/r}(Z)}^{p/r} \Big)^{1/p}.
\end{align*}
Since $p/r > 1$, the operator $\hlmax$ is bounded on $L^{p/r}(Z)$, so we further get
\begin{align*}
  2^{ns} \big\| \sum_{x \in X_n} T(|u|)(x) \chi_{B(x)}\big\|_{L^p(Z)} & \lesssim \Big(\sum_{k \geq n}2^{(n-k)(Q+s-Q/r)(p\land 1)} \big\| U_k \big\|_{L^{p}(Z)}^{p} \Big)^{1/p} \\
& \lesssim \sum_{k \geq n}2^{(n-k)(Q+s-Q/r)(p\land 1)(\frac1p \land 1)} \|U_k\|_{L^p(Z)},
\end{align*}
and as $Q+s-Q/r > 0$, taking the $\ell^q$-norm over $n \in \integer$ yields the boundedness of $T$.

The second part of the statement follows from the first part together with the fact that the composition operator induced by $\Psi$ is bounded on $\ip$, which in turn is an easy consequence of \eqref{eq:seq-norm} and the doubling property of $\mu$.
\end{proof}

\thm{th:trace-commute}
Let $0 < s < \infty$, $Q/(Q+s) < p \leq \infty$ and $0 < q \leq \infty$.

\smallskip
{\rm (i)} If $u \in \ip$, then $\trace u = 0$ $\mu$-almost everywhere.

\smallskip
{\rm (ii)} If $u$ is a sequence on $X$ such that $|du| \in \ip$, then $u - P\trace u \in \ip$ and
\[
  \|u - P \trace u\|_{\ip} \lesssim \|du\|_{\ip}.
\]
In particular, $\trace u = \trace(P\trace u)$ pointwise $\mu$-almost everywhere.

\smallskip
{\rm (iii)} If $f \in \locint(Z)$ and $|d(Pf)| \in \ip$, then $\trace(Pf) = f$ (with convergence in $\locint(Z)$ and pointwise $\mu$-almost everywhere).
\ethm

\begin{proof}
(i) Let $u \in \ip$. Then if $0 < \epsilon < s $ and $B$ is any ball of $Z$ with radius comparable to $2^{-N}$, $N \in \integer$ we get
\begin{align*}
\|u\|_{\ip} & \geq \sup_{k \geq N} 2^{ks} \big\| \sum_{x \in X_k} |u(x)|\chi_{B(x)}\big\|_{L^p(B)}\\
& \geq c(N,s,\epsilon) \bigg( \sum_{k \geq N} 2^{k\epsilon p} \big\| \sum_{x \in X_k} |u(x)|\chi_{B(x)}\big\|_{L^p(B)}^p\bigg)^{1/p} \\
& \approx c(N,s,\epsilon) \bigg( \int_{B} \sum_{|x| \geq N} \big[2^{|x|\epsilon}|u(x)|\big]^p \chi_{B(x)}(\xi)d\mu(\xi)\bigg)^{1/p}.
\end{align*}
Since $B$ has positive and finite $\mu$-measure, this means that for $\mu$-almost all $\xi \in B$, $c_\xi := \sup\{ 2^{|x|\epsilon} |u(x)| \,:\, |x| \geq N,\,B(x)\owns \xi \}$ is finite. For such $\xi$ and $n \geq N$ we then have $|T_n u(\xi)| \lesssim c_\xi 2^{-n\epsilon}$, and the latter quantity converges to zero as $n \to \infty$. This finishes the proof of part (i) since the ball $B$ was arbitrary.

\smallskip
(ii) Let $u$ be as in the statement. Since the limit defining $\trace u$ converges in $\locint(Z)$, we have
\begin{align*}
  |u(x) - P\trace u(x)| & \leq |u(x)- \dashint_{B(x)} T_{|x|}ud\mu| + \sum_{k \geq |x|} \bigg| \dashint_{B(x)} \big(T_{k+1} u - T_k u\big)d\mu \bigg| \\
& \lesssim \sum_{\substack{|y| \geq |x| \\ B(y) \cap B(x) \neq \emptyset}} \frac{\mu\big(B(y)\big)}{\mu\big(B(x)\big)}|du(y)|
\end{align*}
for all $x \in X$, so Proposition \ref{pr:t-operator} yields the desired conclusion.

\smallskip
(iii) This can easily be verified by using the density of Lipschitz functions in $L^1(Z)$. We refer to
\cite[Lemma 4.2]{BS} for details.
\end{proof}

Let us now state the definition of our Besov space. Note that the definition could already have been given earlier; we shall need the tools developed so far in the proof of Theorem \ref{th:seq-identification} below.

\defin{de:besov}
Let $s \in (0,\infty)$, $p \in (Q/(Q+s),\infty]$ and $q \in (0,\infty]$. Then the Besov space $\dbp$ is the vector space of all functions $f \in \locint(Z)$ such that
\beqla{eq:besov-norm}
  \|f\|_{\dbp} := \|d(Pf)\|_{\ip}
\eeq
is finite.
\edefin

\begin{remark}
(i) The vector space $\dbp$ becomes a quasi-normed space after dividing out the functions $f$ such that right-hand side of \eqref{eq:besov-norm} is zero, i.e.~the functions that are constant $\mu$-almost everywhere. We shall frequently abuse notation by writing $\dbp$ for both this quasi-normed space as well as the vector space of functions described above. This remark also applies to other homogeneous function spaces introduced in the sequel.

\smallskip
(ii) Recall that $d(Pf)$ should (for now) be thought of as function on $X$ taking values on a finite-dimensional vector space; see the introduction. The quantity $\|d(Pf)\|_{\ip}$ above thus stands more precisely for $\| \,|d(Pf)|\,\|_{\ip}$. We shall abuse notation in this way in the sequel.

\smallskip
(iii) As mentioned in the introduction, we shall in Section \ref{se:complex-interpolation} below work with an equivalent quasi-norm of $\dbp$ defined in terms of sequence spaces on $E$.
\end{remark}

\begin{remark}
(i) While the space $\ip$ depends on the exact choice of the hyperbolic filling, the space $\dbp$ does not; two admissible choices of $X$ will yield equivalent quasi-norms, with the equivalence constants independent of these two choices. This is not too hard to see using \eqref{eq:seq-norm}, but we are mostly interested in the smoothness range $s \in (0,1)$, and in this case the equivalence is an easy consequence of Proposition \ref{pr:hajlasz} below.

\smallskip
(ii) In view of future applications, it will be useful to note that our methods allow some flexibility in the parameters of the hyperbolic filling. More precisely, the parameters could be chosen so that $(\xi_x)_{x \in X_n}$ is for all $n$ a subset of $Z$ with pairwise distances $\geq c_1 2^{-n}$ for some fixed constant $c_1$ (independent of $n$), that the radii $r_x$ corresponding to the balls $B(x) := B(\xi_x,r_x)$ ($x \in X_n$) are comparable to $2^{-n}$ uniformly in $x$ and $n$, and that the balls $\big(c_2 B(x)\big)_{x \in X_n}$ cover $Z$ for all $n$ where $c_2 \in (0,1)$ is a fixed constant. The results in this paper remain true under these assumptions.
\end{remark}

The results proven so far allow us to prove the following identification of $\dbp$ with a space of sequences on $X$. To explain the notation, fix $x_0 \in X$, write $\izero$ for the space of sequences $u \in \ip$ such that $u(x_0) = 0$, and write $\differencespacei$ for the space of sequences $u\colon X \to \complex$ such that $|du| \in \ip$ and $u(x_0) = 0$.

\thm{th:seq-identification}
Let $0 < s < \infty$, $Q/(Q+s) < p \leq \infty$ and $0 < q \leq \infty$.

\smallskip
{\rm (i)} The trace $\trace u$ of a sequence in $\differencespacei$ is zero $\mu$-almost everywhere if and only if $u \in \izero$.

\smallskip
{\rm (ii)} We have the following isomorphism of quasi-normed spaces:
\[
  \dbp \approx \differencespacei \slash \izero.
\]
In particular, $\dbp$ is a quasi-Banach space, and when $1 < p,\,q < \infty$, it is a reflexive Banach space.
\ethm

\begin{proof}
(i) If $u \in \izero$, then $\trace u = 0$ by Theorem \ref{th:trace-commute} (i). On the other hand, if $u \in \differencespacei$ and $\trace u = 0$, then $u \in \ip$ by Theorem \ref{th:trace-commute} (ii), and since $u(x_0) = 0$, we get $u \in \izero$.

\smallskip
(ii) We first verify that $\differencespacei$ is a quasi-Banach space (in general), and a reflexive Banach space when $1 < p,\,q < \infty$. Using \eqref{eq:seq-norm}, one easily checks that if $(u_k)_{k \geq 1}$ is a Cauchy sequence in $\differencespacei$, then the limit $\lim_{k \to \infty} (u_{k}(x) - u_k(x'))$ exists whenever $x$ and $x'$ are points of $X$ joined by an edge in $E$, and since $u_k(x_0) = 0$ for all $k$, this means that the limit $\lim_{k \to\infty} u_k$ exists pointwise in $X$. It is then not hard to check that this convergence also takes place in $\differencespacei$. Then if $1 < p,\,q < \infty$, it is easily checked using \eqref{eq:seq-norm} that there is a linear bi-Lipschitz embedding from $\differencespace$ to a subspace of a reflexive Banach space of the type $\ell^q(\ell^p)$. By the first part of this proof, this subspace must be closed, and hence reflexive.

Concerning the stated isomorphism, we only note that by Lemma \ref{le:trace}, the operator $\trace\colon \differencespacei \to \locint(Z)$ is continuous, where $\locint(Z)$ is equipped with the topology induced by the seminorms $f \mapsto \|f\|_{B(x_0)}$ and $f \mapsto \inf_{c \in \complex}\|f - c\|_{B(x)}$, $x \in X\backslash\{x_0\}$. Thus $\izero = \trace^{-1}\{0\}$ is a closed subspace of $\differencespacei$, and the quotient space appearing in the statement is well-defined. The actual isomorphism can then be constructed exactly as in the proof of \cite[Theorem 2.9]{BSS}, so we omit the construction here.
\end{proof}

\section{Equivalence of definitions and density of Lipschitz functions}

We shall next show that in the smoothness range $0 < s < 1$, the spaces $\db^s_{p,q}$ coincide with the spaces $\dn^s_{p,q}$ introduced by Koskela, Yang and Zhou \cite{KYZ}. This in particular means that in the Euclidean setting, the spaces $\db^s_{p,q}(\real^d)$ with $0 < s < 1$ and $p > d/(d+s)$ coincide with the standard Fourier-analytically defined Besov spaces. We refer to \cite{GKZ} for a variety of other equivalent characterizations of first-order Besov spaces in the setting of doubling metric measure spaces. The assumptions on our space $(Z,d,\mu)$ and the hyperbolic filling $(X,E)$ are the same as in the previous section.

Let us first recall the definition of the spaces $\dn^s_{p,q}$. For measurable functions $f \colon Z\to\complex$ and $0 < s < \infty$, write $\D^s(f)$ for the collection of all \emph{fractional $s$-Haj\l asz gradients} of $f$, i.e.~sequences $\vec{g} := (g_k)_{k \in \integer}$ of measurable functions $g_k\colon Z \to [0,\infty]$ such that for all integers $k$ and $\xi$, $\eta \in Z$ such that $2^{-k-1} \leq d(\xi,\eta) < 2^{-k}$, we have
\[
  |f(\xi) - f(\eta)| \leq d(\xi,\eta)^s\big(g_k(\xi) + g_k(\eta)\big).
\]
$\dnp$, where $0 < s < \infty$ and $0 < p,\,q \leq \infty$, is then defined as the quasi-normed space of functions (modulo additive constants) $f$ such that
\[
  \|f\|_{\dnp} := \inf_{\vec{g} \in \D^s(f)}  \big\|\big(\|g_k\|_{L^p(Z)}\big)_{k\in\integer}\big\|_{\ell^q}
\] 
is finite.

\prop{pr:hajlasz}
Let $0 < s \leq 1$, $Q/(Q+s) < p \leq \infty$ and $0 < q \leq \infty$. Then
\[
  \dbp = \dnp
\]
with equivalent quasi-norms.
\eprop

The smoothness index $s = 1$ is included here mostly as a curiosity; from \cite[Theorem 4.1]{GKZ} it follows that the spaces $\dn^1_{p,q}$ with $q < \infty$ are often trivial (i.e.~contain only constant functions).

\begin{proof}
Most of the necessary technical details are contained in the proof of \cite[Proposition 3.1]{BSS}, so we will not repeat them here.

To establish the embedding $\dnp \subset \dbp$, we first need to note that the functions in $\dnp$ with $0 < s \leq 1$ and $p > Q/(Q+s)$ are locally integrable. This follows from the Poincar\'e-type inequality \cite[Lemma 2.3]{KYZ}, which is stated in the Euclidean setting, but also holds in setting of doubling metric measure spaces, as noted in the proof of \cite[Theorem 4.1]{KYZ}. Then, if $f \in \dnp$ and $\vec{g} := (g_k)_{k \in \integer} \in \D^s(f)$, the proof of \cite[Proposition 3.1]{BSS} establishes that
\[
  |d(Pf)(x)| \lesssim 2^{-|x|\epsilon'} \sum_{j \geq |x| - \sigma} 2^{-j(s-\epsilon')} \hlmax\big( g_j ^{Q/(Q+\epsilon)}\big)(\xi)^{(Q+\epsilon)/Q}
\]
whenever $\xi \in Z$ and $B(x) \owns \xi$; here $\sigma \geq 0$ is some constant and the parameters $\epsilon$ and $\epsilon'$ are chosen so that $\epsilon < \epsilon' < s$ and $p > Q/(Q+\epsilon)$. One then argues as in the proof of Proposition \ref{pr:t-operator} to obtain
\[
  \|f\|_{\dbp} \lesssim \big\|\big(\|g_k\|_{L^p(Z)}\big)_{k\in\integer}\big\|_{\ell^q},
\]
and taking the infimum over admissible $\vec{g}$ yields the desired embedding.

For the other direction, assume that $f \in \dbp$. The proof of \cite[Proposition 3.1]{BSS} then establishes that the sequence $\vec{g} := (g_k)_{k \in \integer}$ defined by
\[
  g_k(\xi) := 2^{ks}\sum_{|x| \geq k} |d(Pf)(x)| \chi_{B(x)}(\xi)
\]
is a constant times an element of $\D^s(f)$, so
\[
  \|f\|_{\dnp} \lesssim \big\|\big(\|g_k\|_{L^p(Z)}\big)_{k\in\integer}\big\|_{\ell^q} \lesssim \|f\|_{\dbp}.\qedhere
\]
\end{proof}

We shall next examine the density of Lipschitz functions in the spaces $\dbp$. In the previous section, it was shown that the discrete convolutions $T_n(Pf)$ of a function $f \in \dbp$ converge to $f$ in $\locint(Z)$. The following result shows that when $0 < s < 1$ and $q < \infty$, the convergence in fact takes place in $\dbp$.

\thm{th:lipschitz}
{\rm (i)} Let $0 < s < 1$, $Q/(Q+s) < p \leq \infty$ and $0 < q \leq \infty$. Then
\[
  \|T_n(Pf)\|_{\dbp} \lesssim \|f\|_{\dbp}
\]
for all $f \in \dbp$ and $n \in \integer$, with an implied constant independent of $f$ and $n$.

\smallskip
{\rm (ii)} Let $0 < s < 1$, $Q/(Q+s) < p \leq \infty$ and $0 < q < \infty$. Then $T_n(Pf) \to f$ in $\dbp$ as $n \to \infty$.
\ethm

\begin{proof}
We shall use the following two estimates, both established in the proof of \cite[Theorem 3.3]{BSS}:
\beqla{eq:convergence-estimate-1}
  |dP(T_nPf)(x)| \lesssim 2^{n-|x|} \sum_{|y| = n} |d(Pf)(y)|\chi_{B(y)}(\xi) \quad \text{when } |x| \geq n+3 \text{ and } \xi \in B(x)
\eeq
and
\beqla{eq:convergence-estimate-2}
  |dP(f - T_nPf)(x)| \lesssim (T\circ \Psi)(|d(Pf)|\chi_{\cup_{k\geq n}X_k})(x) \quad \text{when } |x| \leq n+2,
\eeq
where $T$ and $\Psi$ are as in Proposition \ref{pr:t-operator}.

\smallskip
(i) From \eqref{eq:convergence-estimate-2} we get $|dP(T_nPf)(x)| \lesssim |d(Pf)(x)| + (T\circ \Psi)(|d(Pf)|\chi_{\cup_{k\geq n}X_k})(x)$ for $|x| \leq n+2$, and we get the desired estimate by combining this with Proposition \ref{pr:t-operator} and \eqref{eq:convergence-estimate-1}.

\smallskip
(ii) From \eqref{eq:convergence-estimate-1} we get $|dP(f - T_nPf)(x)| \lesssim |d(Pf)(x)| + 2^{n-|x|} \sum_{|y| = n} |d(Pf)(y)|\chi_{B(y)}(\xi)$ when $|x| \geq n+3 \text{ and } \xi \in B(x)$, and combining this with \eqref{eq:convergence-estimate-2} and \ref{pr:t-operator} yields
\[
  \|f-T_nPf\|_{\dbp} \lesssim  \bigg(\sum_{k \geq n}2^{ksq}\big\|\sum_{x \in X_k} |d(Pf)(x)| \chi_{B(x)} \big\|_{L^p(Z)}^q\bigg)^{1/q}.
\]
Since $q < \infty$, the latter quantity converges to zero as $n \to \infty$.
\end{proof}

Let us still record the following more abstract density result. Part (i) follows immediately from the Theorem above. Part (ii) requires some additional observations, and its proof will be postponed until after the proof of Proposition \ref{pr:retraction}.

\cor{co:lipschitz} Suppose that $0 < s < 1$ and $0 < q < \infty$.

\smallskip
{\rm (i)} If $Z$ is bounded and $Q/(Q+s) < p \leq \infty$, then Lipschitz functions are dense in $\dbp$.

\smallskip
{\rm (ii)} If $Z$ is unbounded and $Q/(Q+s) < p < \infty$, then Lipschitz functions with bounded support are dense in $\dbp$.
\ecor

In the context of the so-called ``reverse doubling'' metric measure spaces, a result similar to part (ii) above has been obtained in \cite[Proposition 5.21]{HMY}.

\section{Complex interpolation}\label{se:complex-interpolation}

The main aim of this section establish a fairly general complex interpolation formula for the spaces $\dbp$, which appears to be new in the setting of doubling metric measure spaces. For this reason, we will first compute the \emph{Calder\'on product} of two sequence spaces of the type $\ip$; the Calder\'on product is an alternative complex-type interpolation method for quasi-Banach lattices. We recall that in general, the Calder\'on product $X_0^{1-\theta}X_1^\theta$, $0 < \theta < 1$, of two quasi-Banach lattices $X_0$ and $X_1$ on a measure space $M$ is defined as the space of measurable functions $f$ on $M$ such that the quasi-norm
\[
  \|f\|_{X_0^{1-\theta}X_1^{\theta}} := \inf_{ \|g\|_{X_0} \leq 1, \; \|h\|_{X_1} \leq 1} \esssup_{x \in M} \frac{|f(x)|}{|h(x)|^{1-\theta}|g(x)|^\theta}
\]
is finite.

The following result is an analog to \cite[Proposition 6.1]{BSS} concerning sequence spaces related the Triebel-Lizorkin spaces $\dfp$, which in turn was an analog to the Euclidean result \cite[Theorem 8.2]{FJ}. We note that the proof of the Proposition below is much easier than the proofs of \cite[Proposition 6.1]{BSS} and \cite[Theorem 8.2]{FJ} since by \eqref{eq:seq-norm} we can avoid integration on $Z$ completely.

\prop{pr:calderon-product}
Let $0 < s_0,\,s_1 < \infty$, $0 < p_0,\,q_0,\,p_1,\,q_1 \leq \infty$ and $0 < \theta < 1$. Then 
\[
  \big(\I^{s_0}_{p_0,q_0}(X)\big)^{1-\theta} \big(\I^{s_1}_{p_1,q_1}(X)\big)^{\theta} = \ip,
\]
where $s$, $p$ and $q$ are defined by
\beqla{eq:interpolation-params}
      \frac1p = \frac{1-\theta}{p_0} + \frac{\theta}{p_1},
\quad \frac1q = \frac{1-\theta}{q_0} + \frac{\theta}{q_1}
\quad \text{and} \quad s = (1 -\theta)s_0 + \theta s_1.
\eeq
\eprop

\begin{proof}
Suppose first that $u$ belongs to the Calder\'on product space in question. Then if $\lambda$ is any positive number slightly greater than the norm of $u$, we can take $u_0 \in \I^{s_0}_{p_0,q_0}(X)$ and $u_1 \in \I^{s_1}_{p_1,q_1}(X)$ so that $|u(x)| \leq \lambda|u_0(x)|^{1-\theta}|u_1(x)|^{\theta}$ for all $x$ and $\|u_i\|_{\I^{s_i}_{p_i,q_i}}(X) \leq 1$. Using \eqref{eq:seq-norm} and applying H\"older's inequality twice we get
\[
  \|u\|_{\ip} \lesssim \lambda \|u_0\|_{\I^{s_0}_{p_0,q_0}(X)}^{1-\theta} \|u_1\|_{\I^{s_1}_{p_1,q_1}(X)}^\theta \leq \lambda,
\]
and taking the infimum over admissible $\lambda$ shows that the Calder\'on product space is embedded continuously into $\ip$.

For the other direction, we first consider the case with $\min(p_0,p_1) < \infty$ and $\min(q_0,q_1) < \infty$, so that $p,\,q < \infty$. For $u \in \ip\backslash\{0\}$, define the sequences $u_0$, $u_1\colon X\to \complex$ by
\[
  u_i(x) = 2^{\frac{q}{q_i}s|x| - s_i|x|}|u(x)|^{\frac{p}{p_i}} \Big( \sum_{y \in X_{|x|}} \mu\big(B(y) \big)|u(y)|^p \Big)^{\frac{q}{q_ip} - \frac{1}{p_i}},
\]
with the interpretation that $u_i(x) = 0$ if the sum inside the brackets is zero. A simple calculation then shows that $|u(x)| = u_0(x)^{1-\theta} u_1(x)^{\theta}$ for all $x \in X$ and
\[
  \|u_0\|_{\I^{s_0}_{p_0,q_0}(X)}^{1-\theta} \|u_1\|_{\I^{s_1}_{p_1,q_1}(X)}^{\theta} \lesssim \|u\|_{\ip}^{\frac{q}{q_0}(1-\theta)} \|u\|_{\ip}^{\frac{q}{q_1}\theta} = \|u\|_{\ip},
\]
so the desired embedding follows easily.

Suppose now that $\min(p_0,p_1) < \infty$ and $q_0 = q_1 = \infty$. For $u \in \ip\backslash\{0\}$, define $u_0$ and $u_1$ by
\[
  u_i(x) = 2^{s|x| - s_i|x|} |u(x)|^{\frac{p}{p_i}} \Big( \sum_{y \in X_{|x|}} \mu\big(B(y) \big)|u(y)|^p \Big)^{\frac{1}{p} - \frac{1}{p_i}}.
\]
Then $|u(x)| = u_0(x)^{1-\theta}u_1(x)^\theta$ for all $x$ and
\[
  \|u_0\|_{\I^{s_0}_{p_0,q_0}(X)}^{1-\theta} \|u_1\|_{\I^{s_1}_{p_1,q_1}(X)}^{\theta} \lesssim \|u\|_{\ip}^{\frac{p}{p_0}(1-\theta)} \|u\|_{\ip}^{\frac{p}{p_1}\theta} = \|u\|_{\ip},
\]
so the desired embedding again follows.

Finally, the cases with $p_0 = p_1 = \infty$ are relatively easy, so we omit the details.
\end{proof}

It turns out that the Calder\'on product space of two sequence spaces as in the preceding Proposition coincides with the interpolation space $[\I^{s_0}_{p_0,q_0}(X),\I^{s_1}_{p_1,q_1}(X)]_\theta$ obtained by the classical complex interpolation method as long as $\min(p_0,p_1) < \infty$ and $\min(q_0,q_1) < \infty$. It should be noted that although the classical complex interpolation method a priori only makes sense for Banach spaces, it extends to our sequence spaces $\I^s_{p,q}$ since they are \emph{A-convex} or \emph{analytically convex} in the sense of Kalton et al. We refer to the discussion following the proof of \cite[Proposition 6.1]{BSS} and the references therein for details.

With this in mind, we can state our interpolation result for the spaces $\dbp$. This is an analog to \cite[Theorem 6.2]{BSS}. A similar result in the range $1 < p_i,\,q_i < \infty$ in the setting of Ahlfors regular metric measure spaces can be found in \cite{HS}.

\thm{th:interpolation}
Let $0 < s_0,\,s_1 < 1$, $Q/(Q+s_0) < p_0 \leq \infty$, $Q/(Q+s_1) < p_1 \leq \infty$, $0 < q_0,\,q_1 \leq \infty$, $\min(p_0,p_1) < \infty$ and $\min(q_0,q_1) < \infty$. For $0 < \theta < 1$ we then have
\[
  \big[\db^{s_0}_{p_0,q_0}(Z) , \db^{s_1}_{p_1,q_1}(Z) \big]_{\theta} = \dbp
\]
with the parameters $s$, $p$ and $q$ as in \eqref{eq:interpolation-params}
\ethm

This interpolation formula follows readily from Proposition \ref{pr:calderon-product} and the fact (which will be established shortly) that each function space of the type $\db^s_{p,q}$ is obtained as the \emph{retraction} of a sequence space of the type $\I^s_{p,q}$. With these facts in mind, Theorem \ref{th:interpolation} can be proved exactly in the same way as \cite[Theorem 6.2]{BSS}, so we omit the details.

To obtain the function space $\dbp$ as the retraction of a sequence space, it will be convenient for us to work with sequences defined on the edges $E$ of the graph $(X,E)$ rather than on the vertices $X$, since the discrete derivative of a sequence defined on $X$ is naturally defined as a scalar-valued function on $E$. To give the definition of $\ipe$, we first need to introduce some additional notation and conventions.

We equip the edges in $E$ with an orientation, which is chosen so that if $x \sim x'$ and $|x'| > |x|$, then $x'$ is the endpoint of the edge joining $x$ and $x'$. We denote by $e_{x,x'}$ the directed edge from $x$ to $x'$ for all admissible vertices $x$ and $x'$. For a sequence $v\colon X\to \complex$, we write $dv$ for the sequence defined on $E$ by $dv(e_{x,x'}) := v(x') - v(x)$ for all admissible $x$ and $x'$.

For an edge $e \in E$ joining the vertices $x$ and $x'$, write $|e| := |e|\land |e'|$ and $B(e) := B(x)\cup B(x')$. $\ipe$ is then defined as the quasi-normed space of sequences $u \colon E \to \complex$ such that
\begin{align*}
  \|u\|_{\ipe} & := \bigg(\sum_{k\in\integer} 2^{ksq} \big\| \sum_{|e| = k} |u(e)| \chi_{B(e)} \big\|_{L^p(Z)}^q\bigg)^{1/q} \notag \\
& \approx \bigg( \sum_{k\in\integer} 2^{ksq }\Big(\sum_{|e| = k} \mu\big(B(e)\big) |u(e)|^p \Big)^{q/p} \bigg)^{1/q}
\end{align*}
is finite. Proposition \ref{pr:calderon-product} continues to hold with $E$ in place of $X$, as can be easily seen from the proof.

\prop{pr:retraction}
Suppose that $0 < s < 1$, $Q/(Q+s) < p \leq \infty$ and $0 < q \leq \infty$. Then there exist bounded linear operators
\[
  S\colon \dbp \to \ipe \quad \text{and} \quad R\colon \ipe \to \dbp
\]
such that $R\circ S$ is the identity mapping on $\dbp$. More explicitly, if $\xi_0$ is an arbitrary fixed point of $Z$, we may take
\[
  S := f \mapsto d(Pf)
\]
and
\[
  R := u \mapsto \lim_{N \to \infty} \bigg( \sum_{n = -N}^N I_n u(\cdot) - \sum_{n = -N}^{-1} I_n u(\xi_0)  \bigg) + \complex,
\]
where
\[
  I_n u := \sum_{(y,y') \in (X_n \times X_{n+1},\, y\sim y')} u(e_{y,y'}) \psi_{y} \psi_{y'}
\]
for all $n \in \integer$.
\eprop

\begin{proof}
The operator $S$ is bounded between the spaces in question more or less by definition, so the main task is to verify that the operator $R$ is well-defined and bounded, and that $R\circ S$ is the identity operator on $\dbp$. In essence, we have to ``integrate'' sequences defined on $E$ that are not necessarily discrete derivatives of sequences defined on $X$. Most of the required technical work has in fact already been done in the proof of \cite[Proposition 6.3]{BSS}.

Let $u \in \dbp$, and with $\xi_0$ as in the statement of the result, write $\integral u$ for the function defined by
\[
  \integral u(\xi) := \bigg( \sum_{n = -N}^N I_n u(\xi) - \sum_{n = -N}^{-1} I_n u(\xi_0)  \bigg);
\]
that this limit exists in $\locint(Z)$ and pointwise $\mu$-almost everywhere can be verified by using the Lipschitz continuity of the functions $\psi_x$ for the terms with $n < 0$, and by argument similar to the proof of Lemma \ref{le:trace} for the remaining terms (see the proof of \cite[Proposition 6.3]{BSS} for details). Also, write
\[
  U_k := \bigg( \sum_{|e| = k} \big[2^{|e|s}|u(e)|\big]^p \chi_{B(e)}(\cdot) \bigg)^{1/p}
\]
for all $k \in \integer$, so that the $\ipe$-quasi-norm of $u$ is obtained as the mixed $\ell^q(L^p)$-quasi-norm of the functions $U_k$. The proof of \cite[Proposition 6.3]{BSS} then establishes that 
\[
  \big[2^{|x| s} |d(P \integral u)(x)|  \big]^p \lesssim \sum_{n \in \integer} 2^{-\lambda |n-|x||} \hlmax\big( U^r_k\big)(\xi)^{p/r}
\]
for all $\xi \in Z$ and $x \in X$ such that $B(x) \owns \xi$, where $r$ is a parameter strictly between $0$ and $p$, and $\lambda$ is a positive constant. Arguing as in the proof of Proposition \ref{pr:t-operator}, one then gets that
\[
  \| \integral u \|_{\dbp} \lesssim \| u \|_{\ipe}.
\]
Finally, it is easily checked that if $v$ is a sequence on $X$ such that $dv \in \ipe$, then $\integral (dv) = \trace v - T_0 v(\xi_0)$ pointwise $\mu$-almost everywhere, so $R := u \mapsto \integral u + \complex$ satisfies the properties listed in the statement of this Proposition.
\end{proof}

Let us finally show how the second part of Corollary \ref{co:lipschitz} follows from the above Proposition.

\begin{proof}[Proof of Corollary \ref{co:lipschitz} (ii)]
Since $\max(p,q) < \infty$, the vector space $F$ of finitely supported sequences is a dense subspace of $\ipe$. Taking $R$ as in Proposition \ref{pr:retraction}, we then have that $R(F)$ is a dense subspace of $\dbp$, and the elements of $R(F)$ are (equivalence classes of) Lipschitz functions with bounded support.
\end{proof}

\begin{acknowledgements}
The author would like to thank Eero Saksman for reading the manuscript and making several valuable remarks.
\end{acknowledgements}

\end{document}